\documentclass{ifacconf}

\usepackage{mathrsfs}
\usepackage{amsmath}

\usepackage{graphicx}      
\usepackage{natbib}        

\usepackage{balance}

 \usepackage{color} 
\usepackage{amsfonts} 

\usepackage{MnSymbol}%
\usepackage{wasysym}%
\newcommand*{\qedb}{\hfill\ensuremath{\square}}%

\begin{document}
\begin{frontmatter}

\title{A Suboptimality Approach to \\ Distributed $\mathcal{H}_2$ Optimal Control} 

%

\author[First]{Junjie Jiao} 
\author[First]{Harry L. Trentelman} 
\author[First]{M. Kanat Camlibel}

\address[First]{Johann Bernoulli Institute for Mathematics and Computer Science, University of Groningen, Groningen, 9700 AV, The Netherlands (e-mail: \{j.jiao, h.l.trentelman, m.k.camlibel\}@rug.nl).}

\begin{abstract}                
This paper deals with the distributed $\mathcal{H}_2$ optimal control problem for linear multi-agent systems.
In particular, we consider a suboptimal version of the distributed $\mathcal{H}_2$ optimal control problem.
Given a linear multi-agent system with identical agent dynamics and an associated $\mathcal{H}_2$ cost functional, our aim is to design a distributed diffusive static protocol such that the protocol achieves state synchronization for the controlled network and such that the associated cost is smaller than an a priori given upper bound.
We first analyze the $\mathcal{H}_2$ performance of linear systems and then apply the results to linear multi-agent systems.
Two design methods are provided to compute such a suboptimal distributed protocol.
For each method, the expression for the local control gain involves a solution of a single Riccati inequality of dimension equal to the dimension of the individual agent dynamics, and the smallest nonzero and the largest eigenvalue of the graph Laplacian.
\end{abstract}


\end{frontmatter}

\section{Introduction}\label{sec_intro}
The design of distributed protocols for multi-agent systems has received extensive  attention in the past decade (\cite{Olfati-Saber2004}).
This increase in attention is partly due to the broad range of applications of multi-agent systems, e.g. formation control (\cite{Oh2015}), intelligent transportation systems (\cite{Besselink2016}), and smart grids (\cite{Doerfler2013}). 
One of the challenging problems in the context of multi-agent systems is to develop optimal distributed diffusive protocols to minimize given cost performances, while the agents of the network reach a common goal, e.g. state synchronization.
The difficulties of designing such optimal distributed diffusive  protocols are due to the structural constraints on the communication among these agents,  that is, each agent can only receive information from certain other agents.
Therefore, in general, optimal distributed control problems are non-convex and difficult to solve.

To overcome this problem, much effort has been devoted to the design of suboptimal distributed  protocols for multi-agent systems.
In \cite{tamas2008}, the authors established a design method to compute suboptimal distributed controllers subject to a global linear quadratic cost functional. 
The solution of a single LQR problem is required for computing such a suboptimal distributed controller.
Later on, an inverse optimal control problem was addressed in \cite{kristian2014} for both leader-follower and leaderless multi-agent systems.
The authors showed that there exists a global optimal controller if the weighting matrices of the linear quadratic cost functional are chosen to be of a special form.
For other papers related to optimal distributed  control, see also \cite{Mosebach2014}, \cite{nguyen_2017} and \cite{Jiao2018}.

On the other hand, there has been some work on the design of structured controllers for large-scale systems.
In \cite{Rotkowitz2006}, the aim was to design optimal decentralized  controllers, subject to some constraints on the controller structure, to minimize the closed-loop norm of a feedback system.
The authors showed that if the constraints on the controller structure have the property of quadratic invariance, the solution of such problems can be computed efficiently via convex programming.
%
%
%
%
In more recent work, \cite{Fazelnia2017} studied the  distributed optimal problem for linear discrete-time deterministic and stochastic systems. The authors showed that the problem can be relaxed to a semidefinite program, and a globally optimal distributed controller can be obtained if the semidefinite program relaxation has a rank one solution.
In \cite{Fattahi2015}, the authors derived a condition under which, given a optimal centralized controller, there exists a suboptimal distributed controller whose state and input trajectories are close to those of the closed-loop system by using this centralized controller.

In this paper, we study the distributed $\mathcal{H}_2$ optimal control problem for linear multi-agent networks.
We consider a group of identical agents whose dynamics are represented by a finite dimensional linear input/state/output system and a connected, simple undirected weighted graph representing the communication among these agents.
By interconnecting these agents using a distributed diffusive static protocol, we further introduce an $\mathcal{H}_2$ cost functional that penalizes the $\mathcal{L}_2$-norm of the impulse response matrix of the network from the disturbance input to a network output whose components are the weighted differences between the outputs of the individual agent and their neighbors.
The distributed $\mathcal{H}_2$ optimal control problem is then to find the optimal distributed diffusive static protocol that achieves state synchronization for the network and that minimizes the associated cost functional. 
Due to the non-convexity property of the distributed $\mathcal{H}_2$ optimal control problem, this problem is difficult to solve in general.
Therefore, instead of solving the distributed $\mathcal{H}_2$ optimal control problem, we address a suboptimal version of this problem.
More specifically, our aim is to design a distributed diffusive static protocol to achieve state synchronization for the network and to guarantee the associated cost to be smaller than an a priori given upper bound.

The outline of this paper is as follows. 
Section \ref{sec_preliminary} provides some notation and preliminaries on graph theory that will be used throughout this paper. 
In Section \ref{sec_prob}, we formulate the suboptimal distributed $\mathcal{H}_2$ control problem for linear multi-agent systems.
We then present the analysis and design of suboptimal $\mathcal{H}_2$ control for general linear systems in Section \ref{sec_linear_system}, providing necessary results for treating the suboptimal distributed  $\mathcal{H}_2$  control problem.
In Section \ref{sec_state_feed}, we deal with the suboptimal distributed   $\mathcal{H}_2$ control problem for linear multi-agent systems.
Finally, Section \ref{sec_conclusion} concludes this paper.

\section{Preliminaries}\label{sec_preliminary}
%
\subsection{Notation}\label{subsec_notations}
We denote by $\mathbb{R}$ the field of real numbers.
The linear space of real column vectors is denoted by $\mathbb{R}^n$ and the space of real matrices with dimension $m \times n$ is denoted by $\mathbb{R}^{m\times n}$.
Let $\textbf{1}_n\in \mathbb{R}^n$ denote the all-ones vector.
The transpose of a vector $x$ and matrix $X$ are denoted by $x^{\top}$ and $X^{\top}$, respectively.
The inverse of a square matrix $X$ is denoted by $X^{-1}$.
The identity matrix of dimension $n\times n$ is denoted by $I_n$.
For a given symmetric matrix $P$, we write $P>0$ if it is positive definite and $P \geq 0$ if it is positive semidefinite.
The trace of a square matrix $A$ is denoted by $\text{tr}(A)$.
A matrix is called Hurwitz if all its eigenvalues have negative real parts. 
We denote by $\text{diag}(d_1, d_2, \ldots, d_n)$ the $n \times n$ diagonal matrix with $d_1, d_2,\ldots, d_n$ on the diagonal.
%
%
The Kronecker product of two matrices $A \in \mathbb{R}^{m \times n}$ and $B \in \mathbb{R}^{p \times q}$ is denoted by $A \otimes B$ and it has the properties that $(A \otimes B)^{\top} = A^{\top} \otimes B^{\top}$ and $(A_1 \otimes B_1)(A_2 \otimes B_2) = A_1 A_2 \otimes B_1 B_2$ whenever the involved matrix multiplications are legitimate.

\subsection{Graph Theory}\label{subsec_graph}
%
A weighted undirected graph is represented by $\mathcal{G} = (\mathcal{V}, \mathcal{E},\mathcal{A})$, where $\mathcal{V} = \{ 1, 2,\ldots, N \}$ is the node set, $\mathcal{E}$ is the edge set, and $\mathcal{A} = [a_{ij}]$ is the adjacency matrix with nonnegative elements.
The edge set $\mathcal{E}$ of $\mathcal{G}$ is a set of unordered pair $\{ i,j \}$ of distinct nodes  $i$ and $j$ of $\mathcal{G}$, and we have that $a_{ij} > 0$ whenever there is an edge between distinct nodes  $i$ and $j$.
In this paper, we consider simple graphs, i.e. the graphs have no self-loops and hence $a_{ii} = 0$ for all $i$.
%
%
%
%
%
Given a simple weighted undirected  graph $\mathcal{G}$, the degree matrix of $\mathcal{G}$ is the diagonal matrix denoted by $D = \text{diag} ( d_1,d_2,\ldots, d_N )$ with $d_{i} = \sum_{j=1}^{N} a_{ij}$.
Subsequently, we define the Laplacian matrix by $\mathcal{L} = D - \mathcal{A}$.
%
%
The Laplacian matrix $\mathcal{L}$ of an undirected graph is a positive semi-definite symmetric matrix and has real nonnegative eigenvalues. 

A weighted undirected graph is called connected if for each pair of nodes $i$ and $j$ there exists a path from $i$ to $j$.
Furthermore, $\mathcal{G}$ is connected if and only if   $\mathcal{L}$ has a simple eigenvalue $0$. 
In that case, there exists an orthogonal matrix $U$ such that $U^{\top} \mathcal{L} U = \Lambda = \text{diag}(0, \lambda_2, \ldots, \lambda_N)$ with $0 = \lambda_1 < \lambda_2 \leq \cdots \leq \lambda_N$.
For a connected, simple weighted undirected  graph $\mathcal{G}$, let $e_1, e_2,\ldots, e_M$ denote the edges of $\mathcal{G}$, we define the incidence matrix $R \in \mathbb{R}^{N \times M}$ as
\begin{equation*}
R = [r_{ik}], \text{ where } r_{ik} = \left\{
\begin{array}{cl}
1, & \text{ if } e_k =\{ i, j\} \text{ and } i > j,\\
-1, & \text{ if } e_k =\{ i, j\} \text{ and } i < j,\\
0 , & \text{ otherwise},
\end{array} \right.
\end{equation*}
for $i,j = 1,2,\ldots, N$, $i \neq j$ and $k = 1,2, \ldots, M$.
Corresponding to the incidence matrix $R$, we also define the matrix
\begin{equation}\label{W}
	W = \text{diag} ( w_1, w_2,\ldots, w_M)
\end{equation}
as an $M \times M$ diagonal matrix, where $w_i$ is the weight on the edge $i$ for $i = 1,2,\ldots,M$.
The relation between the Laplacian matrix and the incidence matrix is captured by $\mathcal{L} = R W R^{\top}$. See also \cite{Godsil2013}.

\section{Problem Formulation}\label{sec_prob}

In this paper, we consider a multi-agent system consisting of $N$ agents with identical dynamics. 
The interconnection topology among the agents is assumed to be represented by a connected, simple undirected weighted graph with associated graph Laplacian $\mathcal{L}$.
%
The dynamics of agent $i$ is represented by the following continuous-time linear-time-invariant (LTI) system
\begin{equation}\label{mas_decoupled}
\begin{aligned}  
\dot{x}_i (t) &= A x_i (t) + B u_i (t) + E d_i(t), \\
z_i(t) &= C x_i (t) + D u_i(t),
\end{aligned}\quad i = 1,2,\ldots,N
\end{equation}
where $x_i \in \mathbb{R}^n$, $u_i \in \mathbb{R}^m$, $z_i  \in \mathbb{R}^p$ and $d_i \in \mathbb{R}^q$ are the state, the coupling input, the output and the external disturbance input of the $i$th agent, respectively. 
The matrices $A$, $B$, $C$, $D$ and $E$ have suitable dimensions.
%
%
We assume that the pair $(A, B)$ is stabilizable.
%
%
In this paper, we consider the case that the agents (\ref{mas_decoupled}) are interconnected by means of a distributed diffusive static protocol of the form
\begin{equation}\label{protocol_single}
u_i = K \sum_{j = 1}^{N} a_{ij} (x_j - x_i), \quad i =1,2, \ldots, N,
\end{equation}
where $K \in \mathbb{R}^{m \times n}$ is a feedback gain to be designed.

Denote the aggregate vectors as
\begin{align*}
x &= (x_1^{\top}, x_2^{\top}, \ldots,x_N^{\top})^{\top} \in \mathbb{R}^{nN}, 
u = (u_1^{\top}, u_2^{\top}, \ldots, u_N^{\top})^{\top} \in \mathbb{R}^{mN}, \\
z& = (z_1^{\top}, z_2^{\top}, \ldots, z_N^{\top})^{\top} \in \mathbb{R}^{pN}, 
d = (d_1^{\top}, d_2^{\top}, \ldots, d_N^{\top})^{\top} \in \mathbb{R}^{qN}.
\end{align*}
We can then write system (\ref{mas_decoupled}) in compact form as
\begin{equation}\label{mas_compact}
\begin{aligned} 
\dot{x} &= (I_N \otimes A) x + (I_N \otimes B) u + (I_N \otimes E) d ,\\
z &= (I_N \otimes C) x + (I_N \otimes D) u,
\end{aligned}
\end{equation}
the protocol (\ref{protocol_single}) is now of the form
\begin{equation}\label{controller}
u =(\mathcal{L}\otimes K)x.
\end{equation}
Foremost, we want our protocol to achieve state synchronization for the network. This is defined as follows.
\begin{defn}
	The protocol \eqref{controller} is said to achieve state synchronization if, whenever the disturbace input is equal to zero,  i.e. $d=0$, then for all $i = 1, 2, \ldots, N$ we have $x_i(t) - x_j(t) \to 0$ as $t \to \infty$.
\end{defn}
The distributed $\mathcal{H}_2$ optimal control problem is to minimize a given $\mathcal{H}_2$ cost functional for multi-agent system (\ref{mas_compact}) over all protocols (\ref{controller}) that achieve state synchronization. %
Note that in the context of distributed control for multi-agent systems, we are interested in the differences of the state and output values of the agents in the controlled network.
%
%
Observe also that the differences of the state and output values of communicating agents are captured by the incidence matrix $R$ of the underlying graph.
Therefore, we define a new output variable as 
\begin{equation*}
\zeta = (W^{\frac{1}{2}} R^{\top} \otimes I_p )z
\end{equation*}
with $\zeta = (\zeta_1^{\top}, \zeta_2^{\top}, \ldots, \zeta_M^{\top})^{\top} \in \mathbb{R}^{pM}$, where $W$ is the weight matrix given by (\ref{W}). Thus, the output $\zeta$ reflects the weighted disagreement between the outputs of the agents in accordance with the weights of the edges connecting these agents.
Subsequently, we have the following input/state/output model
\begin{equation}\label{network_zeta}
\begin{aligned} 
\dot{x} &= (I_N \otimes A) x + (I_N \otimes B) u + (I_N \otimes E) d ,\\
\zeta &= (W^{\frac{1}{2}}R^{\top}\otimes C)x + (W^{\frac{1}{2}}R^{\top}\otimes D) u.
\end{aligned}
\end{equation}
Next, by substituting protocol (\ref{controller}) into equations (\ref{network_zeta}), we obtain the following equations for the controlled network
\begin{equation*}\label{network_d}
\begin{aligned} 
\dot{x} &= (I_N \otimes A + \mathcal{L} \otimes BK) x + (I_N \otimes E) d , \\
\zeta &= (W^{\frac{1}{2}}R^{\top}\otimes C +  W^{\frac{1}{2}}R^{\top}\mathcal{L} \otimes DK) x.
\end{aligned}
\end{equation*}
Denote
\begin{align*}
&\tilde{A} := I_N \otimes A + \mathcal{L} \otimes BK,\\ &\tilde{E} :=  I_N \otimes E, \\
&\tilde{C} :=W^{\frac{1}{2}}R^{\top}\otimes C +  W^{\frac{1}{2}}R^{\top}\mathcal{L} \otimes DK.
\end{align*}
The impulse response from the disturbance $d$ to the output $\zeta$ is then given by
\begin{equation*}
T_K(t) = \tilde{C}e^{\tilde{A}t}\tilde{E},
\end{equation*}
Subsequently, we define the associated $\mathcal{H}_2$ cost functional as
\begin{equation}\label{cost_K1}
J(K) := \int_{0}^{\infty} \text{tr}\left[ T_K^{\top}(t) T_K(t)\right] dt.
\end{equation}

Recall that our aim is to find a distributed static protocol (\ref{controller}) that minimizes the cost functional (\ref{cost_K1}) over all protocols that achieve state synchronization.
Unfortunately, due to the special form of the protocol which contains the Kronecker product of the feedback gain $K$ and the graph Laplacian $\mathcal{L}$, the distributed $\mathcal{H}_2$ optimal control problem is non-convex and difficult to solve in general. 
Therefore, instead of trying to solve the distributed $\mathcal{H}_2$ optimal control problem itself, we will address a suboptimal version of this problem. 
More specifically, we want to design a state synchronizing, distributed diffusive, static protocol such that the associated cost is smaller than an a priori given upper bound.
More concretely, the problem we want to address is the following:
\begin{prob}\label{prob1}
	Consider multi-agent system (\ref{mas_compact}), with interconnection topology among the agents represented by a connected, simple undirected weighted graph with associated graph Laplacian $\mathcal{L}$, together with cost functional $J(K)$ given by (\ref{cost_K1}). Let $\gamma > 0$ be a given tolerance. Our aim is to design a matrix $K \in \mathbb{R}^{n \times m}$ such that the distributed diffusive static protocol $u = (\mathcal{L} \otimes K)x$ achieves state synchronization and  $J(K) < \gamma$.
	\end{prob}

Before we address Problem \ref{prob1}, we will first briefly discuss the suboptimal $\mathcal{H}_2$  control problem for general linear systems, in that way collecting the required preliminary results to treat the actual suboptimal distributed $\mathcal{H}_2$  control problem for multi-agent systems. 
This will be the subject of the next section.

\section{Suboptimal $\mathcal{H}_2$ Control for Linear Systems} \label{sec_linear_system}

In this section, we consider the suboptimal $\mathcal{H}_2$ control problem for linear systems. 
We will first analyze the $\mathcal{H}_2$ performance of a given system with disturbance inputs. 
Subsequently, we will discuss how to design suboptimal protocols for a linear system with control inputs and disturbance inputs. 
\subsection{$\mathcal{H}_2$ Performance Analysis for Systems with Disturbance Inputs }\label{subsec_d}
In this subsection, we will analyze the $\mathcal{H}_2$ performance for systems with disturbance inputs.
More specifically, we consider the following linear input/state/output system
\begin{equation}\label{sys_d}
\begin{aligned}
\dot{x} (t) &= \bar{A} x(t) + \bar{E} d(t),\\
z(t) &= \bar{C}x (t) 
\end{aligned}
\end{equation}
where $x \in \mathbb{R}^n$ represents the state, $d \in \mathbb{R}^q$ the disturbance input and $z \in \mathbb{R}^p$ the output. 
The matrices $\bar{A}$, $\bar{C}$ and $\bar{E}$ have suitable dimensions.
The impulse response matrix  of system (\ref{sys_d}) from the disturbance $d$ to the output $z$ is $$T(t) = \bar{C} e^{\bar{A}t}\bar{E}.$$
The associated $\mathcal{H}_2$ performance is given by
\begin{equation}\label{cost_d}
J = \int_{0}^{\infty} \text{tr} \left[ T^{\top}(t) T(t) \right] dt,
\end{equation}
which measures the performance of system (\ref{sys_d}) as the square of the $\mathcal{L}_2$-norm of its impulse response matrix.
Note that performance (\ref{cost_d}) is finite if the system is internally stable, i.e., $\bar{A}$ is Hurwitz. 
Our aim is to find conditions such that the performance (\ref{cost_d}) is smaller than a given upper bound.
For this, we have the following lemma. See also \cite{kemin_zhou_book} or \cite{SATO1999295}.
\begin{lem}\label{sys_d_lem}
	Consider system (\ref{sys_d}) with associated performance (\ref{cost_d}).
	The performance is finite if $\bar{A}$ is Hurwitz.
	In that case, we have
	\begin{equation}\label{J_d}
	J = \textnormal{tr}\left(\bar{E}^{\top} Y \bar{E} \right)
	\end{equation}
	where $Y$ is the unique positive semidefinite solution of
	\begin{equation}\label{lyapu_d}
	\bar{A}^{\top} Y + Y \bar{A}  + \bar{C}^{\top} \bar{C} = 0.
	\end{equation} 
	Alternatively, 
	\begin{equation}\label{J_d_ineq}
	J = \inf \{ 
	\textnormal{tr}\left( \bar{E}^{\top} P \bar{E} \right) \ | \  P  \geq 0 
	\text{ and } \bar{A}^{\top} P + P \bar{A}  + \bar{C}^{\top} \bar{C} < 0
	\}.
	\end{equation}
\end{lem}

\begin{pf}
	The fact that the performance (\ref{cost_d}) is given by the expression (\ref{J_d}) involving the Lyapunov equation (\ref{lyapu_d}) is a well-known result.

	Next, we will prove (\ref{J_d_ineq}).
	Let $Y \geq 0$ be the solution of the Lyapunov equation (\ref{lyapu_d}) and let $P\geq 0$ be a solution of the Lyapunov inequality in (\ref{J_d_ineq}). 
	Define $ P := X + Y$.  Then it holds that
	\begin{equation*}
	\bar{A}^{\top} (X + Y) + (X + Y) \bar{A}  + \bar{C}^{\top} \bar{C} < 0.
	\end{equation*}
	Consequently,
	\begin{equation*}
	\bar{A}^{\top} X + X \bar{A} < 0.
	\end{equation*}
	Since $\bar{A}$ is Hurwitz, it follows that $X>0$.
	Thus, we have $P > Y$ and hence $J \leq \textnormal{tr}\left(\bar{E}^{\top} P \bar{E} \right)$ for any $P \geq 0$  satisfying the Lyapunov inequality.

	Next, we will show that for any $\epsilon >0$ there exists $P_{\epsilon} \geq 0$ satisfying the Lyapunov inequality such that $P_{\epsilon} < Y+\epsilon I$, and consequently $\textnormal{tr}\left(\bar{E}^{\top} P_{\epsilon} \bar{E} \right) \leq J +\epsilon \;\textnormal{tr}(\bar{E}^\top \bar{E})$.  
	Indeed, for given $\epsilon$, one can take $P_{\epsilon}$ equal to the unique positive semi-definite solution of
	\begin{equation*}\label{new_lya}
	\bar{A}^{\top} P + P\bar{A}  + \bar{C}^{\top} \bar{C} +\epsilon I_n= 0,
	\end{equation*}
	then $P_{\epsilon} =  \int_{0}^{\infty} e^{\bar{A}^{\top}t} (\bar{C}^{\top} \bar{C} +\epsilon I)e^{\bar{A}t}  \ dt$, 
	so $P_{\epsilon} \downarrow Y$ as $\epsilon \downarrow 0$. This proves our claim.
\qedb
\end{pf}

The following theorem now establishes a {\em necessary} and {\em sufficient} condition (\cite{IWASAKI1994421}), such that the system (\ref{sys_d}) is stable and, for a given upper bound $\gamma > 0$, the  performance (\ref{cost_d}) satisfies $J < \gamma$.
\begin{thm}\label{sys_d_thm}
	Consider system (\ref{sys_d}) with associated performance (\ref{cost_d}).
	Given $\gamma > 0$. 
	Then $\bar{A}$ is Hurwitz and $J < \gamma$ if and only if there exists a positive semidefinite matrix $P$ satisfying
	\begin{align}
	\bar{A}^{\top} P + P \bar{A}  + \bar{C}^{\top} \bar{C} &< 0, \label{lyapu_d_ineq}\\
	\textnormal{tr}\left( \bar{E}^{\top} P \bar{E} \right) &<\gamma. \label{trace_d}
	\end{align}
\end{thm}

\begin{pf}
	($\Leftarrow$) Let $P \geq 0$ satisfy (\ref{lyapu_d_ineq}). 
	Then $\bar{A}^{\top} P + P\bar{A} < 0$. 
	Note also that $P \geq 0$, which implies that $\bar{A}$ is Hurwitz.  
	%
	If $P \geq 0$ also satisfies (\ref{trace_d}), then it follows from Lemma \ref{sys_d_lem} that $J \leq \textnormal{tr}\left( \bar{E}^{\top} P \bar{E} \right) < \gamma$.
	
	($\Rightarrow$)  If $\bar{A}$ is Hurwitz and $J < \gamma$, it follows again from Lemma \ref{sys_d_lem} that there exists a positive semidefinite solution $P$ to (\ref{lyapu_d_ineq}) and (\ref{trace_d})  such that $J \leq \textnormal{tr}\left( \bar{E}^{\top} P \bar{E} \right) < \gamma$.
\qedb
\end{pf}



%
\subsection{Suboptimal $\mathcal{H}_2$ Control for Linear Systems with Control Inputs and Disturbance Inputs}\label{subsec_general}
In this subsection, we will discuss the suboptimal $\mathcal{H}_2$ control problem for linear systems with control inputs and disturbance inputs.
More specifically, we consider the linear input/state/output system
\begin{equation}\label{sys_general_u}
\begin{aligned}
\dot{x}(t) &= A x(t) + B u(t) + E d(t), \\
z(t) &= C x(t) + D u(t),
\end{aligned}
\end{equation}
where $x \in \mathbb{R}^n$ represents the state, $u \in \mathbb{R}^m$ the control input, $z\in \mathbb{R}^p$ the output, and $d\in \mathbb{R}^q$ the disturbance input.
The matrices $A$, $B$, $C$, $D$ and $E$ have suitable dimensions.
We assume that the pair $(A, B)$ is stabilizable.
Using the static state feedback 
\begin{equation}\label{control}
	u = Kx
\end{equation}
yields the closed-loop system
\begin{equation}\label{sys_general}
\begin{aligned}
\dot{x} &= (A + B K) x + E d, \\
z &= (C + DK) x.
\end{aligned}
\end{equation}
We measure the performance of system (\ref{sys_general}) by considering the square of the $\mathcal{L}_2$-norm of its impulse response matrix. 
Therefore, we define the associated $\mathcal{H}_2$ cost functional as
\begin{equation}\label{cost_general}
J(K) = \int_{0}^{\infty} \text{tr} \left[ T_K^{\top}(t) T_K(t) \right] dt
\end{equation}
where 
$$T_K(t) = (C + DK) e^{(A + B K)t} E$$ 
is the closed-loop impulse response matrix of system (\ref{sys_general}) from the disturbance input $d$ to the output $z$.
Let $\gamma > 0$ be a given upper bound for the cost $J(K)$.
We are interested in finding a static state feedback of the form (\ref{control}) such that $A + BK$ is Hurwitz and the associated cost is smaller than the given upper bound $\gamma$, i.e. $J(K) < \gamma$.

The following theorem yields a sufficient condition for the existence of such a static state feedback and how to compute one.


\begin{thm}\label{thm_single_sys}
	Consider system (\ref{sys_general_u}) with associated cost functional (\ref{cost_general}).
	Let $\gamma > 0$.
	Assume that the pair $(A,B)$ is stabilizable.
	Assume that $D^{\top}C = 0$ and $D^{\top} D = I_m$.
	Suppose that there exists  a positive semi-definite matrix $P$ satisfying
	\begin{align}
	A^{\top} P + P A - P B B^{\top} P + C^{\top} C &< 0, \label{are_general}\\
	\textnormal{tr}\left( E^{\top} P E \right) & < \gamma \label{E_cost_general}.
	\end{align}
	Let $K = -B^{\top} P$. Then $A+BK$ is Hurwitz and $J(K) < \gamma$.
\end{thm}
\begin{pf}
	Substituting $K = -B^{\top} P$ into system (\ref{sys_general}) gives us
	\begin{equation*}\label{sys_general_p}
	\begin{aligned}
	\dot{x} &= (A - B B^{\top} P) x + E d, \\
	z &= (C - D B^{\top} P) x.
	\end{aligned}
	\end{equation*}
	Since $D^{\top}C = 0$ and $D^{\top} D = I_m$, inequality (\ref{are_general}) is equivalent to 
	\begin{equation}\label{lyapu_general}
	(A - B B^{\top} P)^{\top} P + P (A - B B^{\top} P) + (C - D B^{\top} P)^{\top}(C - D B^{\top} P) < 0
	\end{equation}
	%
	Since $P \geq 0$ is a solution of (\ref{are_general}), it also satisfies (\ref{lyapu_general}), which implies that $A - B B^{\top} P$ is Hurwitz. 
	Since (\ref{E_cost_general}) also holds, by taking $\bar{A} = A - B B^{\top} P$, $\bar{C} =  C - D B^{\top} P$ and $\bar{E} = E$, it immediately follows from Theorem \ref{sys_d_thm} that $J(K) < \gamma$.
\qedb
\end{pf}

\begin{rem}
	In Theorem \ref{thm_single_sys}, we have assumed that $D^{\top}C = 0$ and $D^{\top} D = I_m$, which is often called the {\em standard form}. Although we do not consider the general case here, it is straightforward to extend our result to the general case, since the general problem can be reduced to a problem in standard form by a preliminary state feedback transformation. See e.g. \cite{harry_book}.
\end{rem}

%
%
\section{Distributed Suboptimal $\mathcal{H}_2$ Control for Multi-Agent Systems}\label{sec_state_feed}
In the previous section, we have dealt with the suboptimal $\mathcal{H}_2$ control problem for linear systems, collecting the necessary results for treating the suboptimal distributed  $\mathcal{H}_2$ control problem.
In the present section, we deal with the  suboptimal distributed  $\mathcal{H}_2$ control problem for multi-agent networks with identical linear agent dynamics.

As has already been shown in Section \ref{sec_prob}, the input/state/ output model of the multi-agent network we consider is given by 
\begin{equation}\label{network_dd}
\begin{aligned} 
\dot{x} &= (I_N \otimes A + \mathcal{L} \otimes BK) x + (I_N \otimes E) d , \\
\zeta &= (W^{\frac{1}{2}}R^{\top}\otimes C +  W^{\frac{1}{2}}R^{\top}\mathcal{L} \otimes DK) x.
\end{aligned}
\end{equation}
For convenience, we also repeat here the associated $\mathcal{H}_2$ cost functional 
\begin{equation}\label{cost_K}
J(K) = \int_{0}^{\infty} \text{tr}\left[ T_K^{\top}(t) T_K(t)\right] dt,
\end{equation}
where 
$T_K(t) = \tilde{C}e^{\tilde{A}t}\tilde{E}$
is the impulse response matrix from the disturbance input $d$ to the output $\zeta$ with $\tilde{A} := I_N \otimes A + \mathcal{L} \otimes BK$, $\tilde{E} :=  I_N \otimes E$ and
$\tilde{C} :=W^{\frac{1}{2}}R^{\top}\otimes C +  W^{\frac{1}{2}}R^{\top}\mathcal{L} \otimes DK$.

The suboptimal distributed  $\mathcal{H}_2$ control problem is to find a distributed diffusive static protocol (\ref{controller}) with gain matrix $K$ that achieves state synchronization and such that the associated cost (\ref{cost_K}) is smaller than a given upper bound $\gamma >0$, i.e. $J(K) <\gamma$.
We further assume that $D^{\top}C = 0$ and $D^{\top} D = I_m$, i.e. we assume that  the suboptimal distributed  $\mathcal{H}_2$ control problem is in standard form.

We first apply the state transformation 
\begin{equation*}\label{state_trans}
	\bar{x} = (U^{\top} \otimes I_n) x
\end{equation*}
where the orthogonal matrix $U$ is defined in Section \ref{subsec_graph}. 
After this state transformation, the equations of the controlled network become
\begin{align*}
\dot{\bar{x}} &= (I_N \otimes A + \Lambda \otimes BK) \bar{x} + (U^{\top} \otimes E) d, \\
\zeta & = (W^{\frac{1}{2}}R^{\top} U\otimes C +  W^{\frac{1}{2}}R^{\top} \mathcal{L} U \otimes DK) \bar{x}, 
\end{align*}
and our cost functional is equal to
\begin{equation}\label{cost_KK}
	J(K) = \int_{0}^{\infty} \text{tr}\left[ \bar{T}_K^{\top}(t) \bar{T}_K(t)\right] dt,
\end{equation}
where
\begin{equation}\label{inpul_resp_trans}
	\bar{T}_K(t) = C_o e^{A_o t} E_o
\end{equation}
is the impulse response matrix from the disturbance input $d$ to the output $\zeta$ with $A_o := I_N \otimes A + \Lambda \otimes BK$, $C_o := W^{\frac{1}{2}}R^{\top} U\otimes C +  W^{\frac{1}{2}}R^{\top} \mathcal{L} U \otimes DK$ and $E_o := U^{\top} \otimes E$.
Note that, by applying the state transformation, only the system model has changed while the associated cost remains the same.

In order to proceed, we introduce the following input/state/ output systems
\begin{equation}\label{local_sys}
\begin{aligned}
\dot{\xi}_i &= A\xi_i + \lambda_i B v_i + E \delta_i, \\
\eta_i &= \sqrt{\lambda_i} C \xi_i + \lambda_i \sqrt{\lambda_i} D v_i,
\end{aligned}\quad i = 2, 3,\ldots, N.
\end{equation}
where $\lambda_i, i=2,3,\ldots, N$ are the nonzero eigenvalues of the graph Laplacian $\mathcal{L}$.
Using in all systems (\ref{local_sys}) the identical static state feedback
\begin{equation}\label{state_feedback}
	v_i = K \xi_i, \quad i=2,3,\ldots,N
\end{equation}
yields the closed-loop systems
\begin{equation}\label{local_sys_closed}
\begin{aligned}
	\dot{\xi}_i &= (A + \lambda_i B K)\xi_i + E \delta_i, \\
	\eta_i &= (\sqrt{\lambda_i} C + \lambda_i\sqrt{\lambda_i} DK ) \xi_i,
\end{aligned}\quad i = 2, 3,\ldots, N.
\end{equation}
We further introduce the associated cost functionals
\begin{equation}\label{local_cost}
	J_i(K) = \int_{0}^{\infty} \textnormal{tr} \left[ T_{i,K}^{\top}(t) T_{i,K}(t) \right] dt, \quad i = 2, 3,\ldots, N,
\end{equation}
where 	
\begin{equation}\label{tf}
	T_{i,K} = (\sqrt{\lambda_i} C + \lambda_i\sqrt{\lambda_i} D K) e^{(A + \lambda_i B K)t} E, \quad i = 2, 3,\ldots, N
\end{equation}
are the closed-loop impulse response matrices from the disturbance $\delta_i$ to the output $\eta_i$, for $i = 2, 3,\ldots, N$, respectively.

It turns out that our original cost functional can be expressed as the sum of the cost functionals associated with the auxiliary systems (\ref{local_sys}).
In fact, the following theorem holds.

\begin{thm}\label{thm_costcost}
	Consider the network (\ref{network_dd}) with associated cost (\ref{cost_K}) and the systems (\ref{local_sys_closed}) with associated costs (\ref{local_cost}) for $i =2,3,\ldots,N$, respectively.
	Then the protocol (\ref{controller}) achieves state synchronization for the network (\ref{network_dd}) if and only if the static state feedback (\ref{state_feedback}) internally stabilizes all systems (\ref{local_sys}).
	Moreover,
	\begin{equation}\label{cost_equivalence}
		J(K) = \sum_{i=2}^{N} J_i(K).
	\end{equation}
\end{thm}

\begin{pf}
It is a standard result that the protocol (\ref{controller}) achieves state synchronization for the network (\ref{network_dd}) if and only if the static state feedback (\ref{state_feedback}) internally stabilizes all systems (\ref{local_sys}). See e.g. \cite{zhongkui_li_unified_2010} or \cite{harry_2013}. 

We now prove (\ref{cost_equivalence}).
First, note that the cost (\ref{cost_K}) for the network (\ref{network_dd}) is equal to (\ref{cost_KK}).
Then, substituting (\ref{inpul_resp_trans}) into (\ref{cost_KK}) gives us
\begin{align*}
	J(K) &= \int_{0}^{\infty}\text{tr} \left[ (U\otimes E^{\top}  ) 
	e^{(I_N \otimes A + \Lambda \otimes BK)^{\top} t} \right.\nonumber\\ 
	&\qquad (W^{\frac{1}{2}}R^{\top} U\otimes C +  W^{\frac{1}{2}}R^{\top} \mathcal{L} U \otimes DK)^{\top} \nonumber\\
	&\qquad (W^{\frac{1}{2}}R^{\top} U\otimes C +  W^{\frac{1}{2}}R^{\top} \mathcal{L} U \otimes DK)\nonumber \\
	&\qquad \left. e^{(I_N \otimes A + \Lambda \otimes BK) t}(U^{\top} \otimes E) \right] dt,
\end{align*}
which is equal to
\begin{align}
	J(K) &=\int_{0}^{\infty}\text{tr}  \left[ (I_N \otimes E^{\top}) 
	e^{(I_N \otimes A + \Lambda \otimes BK)^{\top} t} \right.\nonumber\\ 
	&\qquad (W^{\frac{1}{2}}R^{\top} U \otimes C +  W^{\frac{1}{2}}R^{\top}\mathcal{L} U \otimes DK)^{\top} \nonumber\\
	&\qquad (W^{\frac{1}{2}}R^{\top} U  \otimes C +  W^{\frac{1}{2}}R^{\top}\mathcal{L} U \otimes DK) \nonumber\\
	&\qquad \left. e^{(I_N \otimes A + \Lambda \otimes BK) t}(I_N \otimes E ) \right] dt.\label{JK}
\end{align}
Recall that $U^{\top} \mathcal{L} U = \Lambda = \text{diag}(0, \lambda_2, \ldots, \lambda_N)$, $\mathcal{L} = R W R^{\top}$ and $D^{\top}C = 0$. 
Therefore, (\ref{JK}) is equal to 
\begin{align*}
J(K) &= \int_{0}^{\infty}\text{tr} \left[ (I_N \otimes E^{\top}) 
	e^{(I_N \otimes A + \Lambda \otimes BK)^{\top} t} \right. \nonumber\\ 
	&\qquad (\Lambda \otimes C^{\top} C +  \Lambda^3 \otimes K^{\top}D^{\top} D K) \nonumber\\
	&\qquad \left. e^{(I_N \otimes A + \Lambda \otimes BK) t}(I_N \otimes E ) \right] dt \nonumber \\
	& = \sum_{i=2}^{N} \int_{0}^{\infty} \text{tr} \left[E^{\top}
	e^{(A + \lambda_i BK)^{\top} t} \right.\nonumber\\ 
	&\qquad (\lambda_i C^{\top} C + \lambda_i^3 K^{\top} D^{\top} D K) \left. e^{(A + \lambda_i BK) t}E \right] dt \nonumber \\
	&=\sum_{i=2}^{N}\int_{0}^{\infty}  \text{tr}  \left[E^{\top}
	e^{(A + \lambda_i BK)^{\top} t} \right.\nonumber\\ 
	&\qquad (\sqrt{\lambda_i} C + \lambda_i\sqrt{\lambda_i} D K)^{\top} \\
	&\qquad (\sqrt{\lambda_i} C +  \lambda_i\sqrt{\lambda_i} D K)\left. e^{(A + \lambda_i BK) t}E \right] dt \nonumber\\
	&= \sum_{i=2}^{N}J_i(K) \nonumber
\end{align*}
\qedb
\end{pf}
%
Based on Theorem \ref{thm_costcost}, we have transformed the problem of suboptimal distributed  $\mathcal{H}_2$ control for the network (\ref{network_dd}) into suboptimal $\mathcal{H}_2$ control problems for $N-1$ linear systems (\ref{local_sys_closed}) with the same feedback gain $K$.
%
%
Next, we want to establish conditions under which all $N-1$ systems (\ref{local_sys_closed}) are internally stable and the state feedback \eqref{state_feedback} is suboptimal.

The following lemma yields a {\em necessary} and {\em sufficient} condition for a given gain matrix $K \in \mathbb{R}^{m \times n}$ such that all systems (\ref{local_sys_closed}) are internally stable and $\sum_{i=2}^{N}J_i(K) <\gamma$.

\begin{lem}\label{lem_mas}
	Consider the closed-loop systems (\ref{local_sys_closed}) with the associated cost functionals \eqref{local_cost}.
	The static state feedback \eqref{state_feedback} internally stabilizes all systems and $\sum_{i=2}^{N}J_i(K) <\gamma$ if and only if there exist positive semidefinite matrices $P_i, i=2,3,\ldots,N$ satisfying
	\begin{align}
	(A + \lambda_i B K)^{\top} P_i + P_i (A + \lambda_i B K) & \nonumber\\
	+ (\sqrt{\lambda_i}  C + \lambda_i\sqrt{\lambda_i}  DK)^{\top}  (\sqrt{\lambda_i}  C& + \lambda_i\sqrt{\lambda_i}  DK)  < 0, \label{lyapu_mas}\\
	\sum_{i=2}^{N} \textnormal{tr}
	& \left( E^{\top} P_i E \right)  < \gamma. \label{E_P}
	\end{align}
\end{lem}

\begin{pf}
	($\Leftarrow$)
	Since (\ref{E_P}) holds, there exist sufficiently small $\epsilon_i > 0, i=2,\ldots,N$ such that $\sum_{i=2}^{N} \gamma_i < \gamma$  where $\gamma_i := \textnormal{tr} \left( E^{\top} P_i E \right) +\epsilon_i$. 
	Because there exists $P_i$ such that (\ref{lyapu_mas}) and $\textnormal{tr} \left( E^{\top} P_i E \right) < \gamma_i$ hold for all $i =2,\ldots,N$, by taking $\bar{A} = A + \lambda_i BK$ and $\bar{C} = \sqrt{\lambda_i}  C + \lambda_i\sqrt{\lambda_i}  DK$, it follows from Theorem \ref{sys_d_thm} that all systems (\ref{local_sys_closed}) are internally stable and $J_i(K) <\gamma_i$ for $i = 2,\ldots,N$. 
	Therefore, $\sum_{i=2}^{N}J_i(K)  < \gamma$.
	
	($\Rightarrow$)
	Since $\sum_{i=2}^{N}J_i(K) <\gamma$, there exist sufficiently small $\epsilon_i > 0, i=2,\ldots,N$ such that $\sum_{i=2}^{N} \gamma_i < \gamma$  where $\gamma_i := J_i(K)+\epsilon_i$.
	Because all systems (\ref{local_sys_closed}) are internally stable and $J_i(K) <\gamma_i$ for $i=2,\ldots, N$, by taking $\bar{A} = A + \lambda_i BK$ and $\bar{C} =\sqrt{\lambda_i}  C + \lambda_i\sqrt{\lambda_i}  DK$, it follows from Theorem \ref{sys_d_thm} that there exist positive semi-definite matrices $P_i$ such that (\ref{lyapu_mas})  and  $\textnormal{tr} \left( E^{\top} P_i E \right) < \gamma_i$ hold for all $i =2,\ldots,N$. 
	Since $\sum_{i=2}^{N} \gamma_i < \gamma$, this implies that $\sum_{i=2}^{N} \textnormal{tr} \left( E^{\top} P_i E \right) < \gamma$.
\qedb
\end{pf}

Lemma \ref{lem_mas} establishes a necessary and sufficient condition for a given gain matrix $K$ to internally stabilize all closed-loop systems (\ref{local_sys_closed}) and to achieve $\sum_{i=2}^{N}J_i(K)<\gamma$. 
However, it does yet not provide a method to compute such gain matrix $K$.
To this end, the following two theorems provide two design methods for computing such a  gain matrix $K$ and, correspondingly, two suboptimal distributed protocols for multi-agent system (\ref{mas_decoupled}) together with cost functional (\ref{cost_K}). 

\begin{thm}\label{main1}
	Consider multi-agent system (\ref{mas_decoupled}) with the associated cost functional (\ref{cost_K}).
	Choose $c$ such that 
	\begin{equation}\label{c1}
	0 < c \leq \frac{2}{\lambda_2^2 + \lambda_2 \lambda_N + \lambda_N^2}.
	\end{equation}
	Then there exists a positive semidefinite matrix $P$ satisfying
	\begin{equation}\label{one_are}
	A^{\top} P + P A + (c^2 \lambda_2^3 - 2 c \lambda_2) PB B^{\top} P +\lambda_N C^{\top} C <0.
	\end{equation}
	Assume, moreover, that $P$ also satisfies
	\begin{equation}\label{p}
		\textnormal{tr} \left( E^{\top} P E \right)  < \frac{\gamma}{N-1}.
	\end{equation} 
	Let $K := -cB^{\top}P$. Then protocol (\ref{controller}) achieves synchronization, and the protocol is suboptimal, i.e. $J(K) < \gamma$.
\end{thm}
\begin{pf}
	Using the upper and lower bound on $c$ given by (\ref{c1}), it can be verified that $c^2 \lambda_i^3 - 2 c \lambda_i \leq c^2 \lambda_2^3 - 2 c \lambda_2 <0$ for  $i=2,3,\ldots,N$. 
	Since also $\lambda_i \leq \lambda_N$, one can see that the positive semidefinite solution $P$ of (\ref{one_are}) also satisfies the $N-1$ Riccati inequalities
	\begin{equation}\label{n_are1}
	A^{\top} P + P A + (c^2 \lambda_i^3 - 2 c \lambda_i) PB B^{\top} P +\lambda_i C^{\top} C <0, \quad i=2,\ldots,N.
	\end{equation}
	Equivalently, $P$ also satisfies the Lyapunov inequalities  
	\begin{equation*}\label{lya_ineq}
	\begin{aligned}
		(A-c\lambda_i B B^{\top}P)^{\top}P & + P(A- c\lambda_iB B^{\top}P) \\
		& + c^2\lambda_i^3 PB B^{\top}P + \lambda_i C^{\top} C < 0,
	\end{aligned}
	\end{equation*}
	for $i=2,\ldots,N$.
	Taking $P_i = P$ for $i = 2,3,\ldots,N$ and $K := -c B^{\top}P$ in inequalities (\ref{lyapu_mas}) and (\ref{E_P}) immediately gives us inequalities (\ref{n_are1}) and 
	\begin{equation}\label{costN}
	\sum_{i=2}^{N} \textnormal{tr} \left( E^{\top} P  E \right)  < \gamma.
	\end{equation}
	%
	Then it follows from Lemma \ref{lem_mas} that all systems (\ref{local_sys_closed}) are internally stable and $\sum_{i=2}^{N}J_i(K)<\gamma$. 
	Furthermore, it follows from  Theorem \ref{thm_costcost} that the protocol (\ref{controller}) achieves state synchronization in the network (\ref{network_dd}) and $J(K)<\gamma$.
\qedb
\end{pf}

\begin{thm}\label{main2}
	Consider the multi-agent system (\ref{mas_decoupled}) with associated cost functional (\ref{cost_K}).
	Choose $c$ such that 
	\begin{equation}\label{c2}
	\frac{2}{\lambda_2^2 + \lambda_2 \lambda_N + \lambda_N^2} < c <\frac{2}{\lambda_N^2}.
	\end{equation}
	Then there exists a positive semidefinite matrix $P$ satisfying
	\begin{equation}\label{one_are1}
	A^{\top} P + P A + (c^2 \lambda_N^3 - 2 c \lambda_N) PB B^{\top} P +\lambda_N C^{\top} C <0.
	\end{equation}
	Assume, moreover, that $P$ also satisfies
	\begin{equation}
	\textnormal{tr} \left( E^{\top} P E \right)  < \frac{\gamma}{N-1}.
	\end{equation} 
	Let $K := -cB^{\top}P$. Then protocol (\ref{controller}) achieves state synchronization, and the protocol is suboptimal, i.e. $J(K) < \gamma$.
\end{thm}
\begin{pf}
	The proof is similar to the proof of Theorem \ref{main1} and hence is omitted here.
\qedb
\end{pf}

\balance

\section{Conclusion}\label{sec_conclusion}
In this paper we have studied a suboptimal distributed  $\mathcal{H}_2$ control problem for linear multi-agent systems with connected, simple undirected weighted graph.
Given a multi-agent system  with identical agent dynamics and an associated global $\mathcal{H}_2$ cost functional, we provide two design methods for computing a suboptimal distributed  static protocol such that the protocol achieves state synchronization for the controlled network  and the associated cost is smaller than a given upper bound.
For each method, the expression for the local control gain is provided in terms of solutions of a single Riccati inequality, whose dimension is equal to the dimension of the individual agent dynamics, and also involves the largest and the smallest nonzero eigenvalue of the graph Laplacian.


\bibliography{h2_ifac}             

\begin{thebibliography}{19}
\providecommand{\natexlab}[1]{#1}
\providecommand{\url}[1]{\texttt{#1}}
\providecommand{\urlprefix}{URL }
\expandafter\ifx\csname urlstyle\endcsname\relax
  \providecommand{\doi}[1]{doi:\discretionary{}{}{}#1}\else
  \providecommand{\doi}{doi:\discretionary{}{}{}\begingroup
  \urlstyle{rm}\Url}\fi

\bibitem[{Besselink et~al.(2016)Besselink, Turri, van~de Hoef, Liang, Alam,
  M{\aa}rtensson, and Johansson}]{Besselink2016}
Besselink, B., Turri, V., van~de Hoef, S.H., Liang, K.Y., Alam, A.,
  M{\aa}rtensson, J., and Johansson, K.H. (2016).
\newblock Cyber--physical control of road freight transport.
\newblock \emph{Proceedings of the IEEE}, 104(5), 1128--1141.

\bibitem[{Borrelli and Keviczky(2008)}]{tamas2008}
Borrelli, F. and Keviczky, T. (2008).
\newblock Distributed {LQR} design for identical dynamically decoupled systems.
\newblock \emph{IEEE Transactions on Automatic Control}, 53(8), 1901--1912.

\bibitem[{D{\"o}rfler et~al.(2013)D{\"o}rfler, Chertkov, and
  Bullo}]{Doerfler2013}
D{\"o}rfler, F., Chertkov, M., and Bullo, F. (2013).
\newblock Synchronization in complex oscillator networks and smart grids.
\newblock \emph{Proceedings of the National Academy of Sciences}, 110(6),
  2005--2010.

\bibitem[{Fattahi et~al.(2015)Fattahi, Fazelnia, and Lavaei}]{Fattahi2015}
Fattahi, S., Fazelnia, G., and Lavaei, J. (2015).
\newblock Transformation of optimal centralized controllers into near-global
  static distributed controllers.
\newblock In \emph{Decision and Control (CDC), 2015 IEEE 54th Annual Conference
  on}, 4915--4922. IEEE.

\bibitem[{Fazelnia et~al.(2017)Fazelnia, Madani, Kalbat, and
  Lavaei}]{Fazelnia2017}
Fazelnia, G., Madani, R., Kalbat, A., and Lavaei, J. (2017).
\newblock Convex relaxation for optimal distributed control problems.
\newblock \emph{IEEE Transactions on Automatic Control}, 62(1), 206--221.

\bibitem[{Godsil and Royle(2013)}]{Godsil2013}
Godsil, C. and Royle, G. (2013).
\newblock \emph{Algebraic Graph Theory}.
\newblock Graduate Texts in Mathematics. Springer New York.

\bibitem[{Iwasaki et~al.(1994)Iwasaki, Skelton, and Geromel}]{IWASAKI1994421}
Iwasaki, T., Skelton, R., and Geromel, J. (1994).
\newblock Linear quadratic suboptimal control with static output feedback.
\newblock \emph{Systems \& Control Letters}, 23(6), 421 -- 430.

\bibitem[{Jiao et~al.(2018)Jiao, Trentelman, and Camlibel}]{Jiao2018}
Jiao, J., Trentelman, H.L., and Camlibel, M.K. (2018).
\newblock A suboptimality approach to distributed linear quadratic optimal
  control.
\newblock \emph{manuscript 2018}, submitted for publication.

\bibitem[{Li et~al.(2010)Li, Duan, Chen, and Huang}]{zhongkui_li_unified_2010}
Li, Z., Duan, Z., Chen, G., and Huang, L. (2010).
\newblock Consensus of multiagent systems and synchronization of complex
  networks: A unified viewpoint.
\newblock \emph{IEEE Transactions on Circuits and Systems I: Regular Papers},
  57(1), 213--224.

\bibitem[{Mosebach and Lunze(2014)}]{Mosebach2014}
Mosebach, A. and Lunze, J. (2014).
\newblock Synchronization of autonomous agents by an optimal networked
  controller.
\newblock In \emph{Proc. European Control Conf. (ECC)}, 208--213.

\bibitem[{Movric and Lewis(2014)}]{kristian2014}
Movric, K.H. and Lewis, F.L. (2014).
\newblock Cooperative optimal control for multi-agent systems on directed graph
  topologies.
\newblock \emph{IEEE Transactions on Automatic Control}, 59(3), 769--774.

\bibitem[{Nguyen(2017)}]{nguyen_2017}
Nguyen, D.H. (2017).
\newblock Reduced-order distributed consensus controller design via edge
  dynamics.
\newblock \emph{IEEE Transactions on Automatic Control}, 62(1), 475--480.

\bibitem[{Oh et~al.(2015)Oh, Park, and Ahn}]{Oh2015}
Oh, K.K., Park, M.C., and Ahn, H.S. (2015).
\newblock A survey of multi-agent formation control.
\newblock \emph{Automatica}, 53, 424--440.

\bibitem[{Olfati-Saber and Murray(2004)}]{Olfati-Saber2004}
Olfati-Saber, R. and Murray, R.M. (2004).
\newblock Consensus problems in networks of agents with switching topology and
  time-delays.
\newblock \emph{IEEE Transactions on Automatic Control}, 49(9), 1520--1533.

\bibitem[{Rotkowitz and Lall(2006)}]{Rotkowitz2006}
Rotkowitz, M. and Lall, S. (2006).
\newblock A characterization of convex problems in decentralized control.
\newblock \emph{IEEE Transactions on Automatic Control}, 51(2), 274--286.

\bibitem[{Sato and Liu(1999)}]{SATO1999295}
Sato, T. and Liu, K.Z. (1999).
\newblock {LMI} solution to general {$H_2$} suboptimal control problems.
\newblock \emph{Systems \& Control Letters}, 36(4), 295 -- 305.

\bibitem[{Trentelman et~al.(2013)Trentelman, Takaba, and
  Monshizadeh}]{harry_2013}
Trentelman, H.L., Takaba, K., and Monshizadeh, N. (2013).
\newblock Robust synchronization of uncertain linear multi-agent systems.
\newblock \emph{IEEE Transactions on Automatic Control}, 58(6), 1511--1523.

\bibitem[{Trentelman et~al.(2001)Trentelman, Stoorvogel, and
  Hautus}]{harry_book}
Trentelman, H.L., Stoorvogel, A.A., and Hautus, M. (2001).
\newblock \emph{Control Theory for Linear Systems}.
\newblock Springer Verlag.

\bibitem[{Zhou et~al.(1996)Zhou, Doyle, and Glover}]{kemin_zhou_book}
Zhou, K., Doyle, J.C., and Glover, K. (1996).
\newblock \emph{Robust and Optimal Control}.
\newblock New Jersey, NJ, USA: Prentice-Hall.

\end{thebibliography}
                                                   







\end{document}